\documentclass[10 pt,twosided,reqno]{amsart}
\usepackage{amscd,amssymb}
\usepackage[arrow,matrix]{xy}
\usepackage{graphicx}
\theoremstyle{plain}
\newtheorem{theorem}[subsection]{Theorem}
\newtheorem{lemma}[subsection]{Lemma}
\newtheorem{proposition}[subsection]{Proposition}
\newtheorem{cor}[subsection]{Corollary}

\theoremstyle{definition}
\newtheorem{rk}[subsection]{Remark}
\newtheorem{defn}[subsection]{Definition}
\newtheorem{ex}[subsection]{Example}

\numberwithin{equation}{section} 

\begin{document}
\title [FIBONACCI MODULES AND MULTIPLE FIBONACCI SEQUENCES]
{FIBONACCI MODULES AND MULTIPLE FIBONACCI SEQUENCES }
\thanks{This research is partially supported by Higher Education Commission, Pakistan.}



\maketitle
 \pagestyle{myheadings} \markboth{\centerline {\scriptsize
A. R. NIZAMI }} {\centerline {\scriptsize FIBONACCI MODULES AND
MULTIPLE FIBONACCI SEQUENCES }}

\begin{center}
{\centerline {  ABDUL RAUF NIZAMI}} \vspace{0.15cm} {\small
\noindent
\emph{Abdus Salam School of Mathematical Sciences,\\
 GC University, Lahore-Pakistan.}\\
       \emph{ e-mail: arnizami@yahoo.com}}
\end{center}
\vspace{0.2cm}
\begin{abstract}
Double Fibonacci sequences $(x_{n,k})$ are introduced and they are
related to operations with Fibonacci modules. Generalizations and
examples  are also discussed.
\end{abstract}


\section{\large{\textbf{Introduction}}} \label{sec1}
Let us fix a commutative ring $\mathcal{R}$; $\mathcal{R}^{2}$ will
denote the rank $2$ free $\mathcal{R}$-module and also the product
ring $\mathcal{R\times R}$. The main object of study is the
Fibonacci module of type $(a,b)\in\mathcal{R}^{2}$ associated to the
$\mathcal{R}$-module $\mathbf{M}$:
\begin{defn}\label{def1.1}
$\mathcal{F}_{\mathbf{M}}(a,b)$ is the set of sequences
$\{(x_{n})_{n\geq0}: x_{n}\in\mathbf{M},\,\,x_{n+2}=ax_{n+1}+b x_{n},\forall n\geq0\}$. If
$\mathbf{M}=\mathcal{R}$, we use the shorter notation
$\mathcal{F}(a,b)$.
\end{defn}

\begin{rk}\label{rem1.2}
Using the $\mathcal{R}[T]$ structure of the $\mathcal{R}$-module of
all sequences in $\mathbf{M}$:
$\mathcal{S}_{\mathbf{M}}=\{(x_{n})_{n\geq0}: x_{n}\in\mathbf{M}\}$,
where the action $T$ is given by the shift
$T(x_{0},x_{1},x_{2},\ldots)=(x_{1},x_{2},x_{3},\ldots)$, one can
describe $\mathcal{F}_{\mathbf{M}}(a,b)$ as the sub
$\mathcal{R}[T]$-module $\ker(T^{2}-aT-b)$. We also consider
$\mathcal{\widetilde{F}}_{\mathbf{M}}(a,b)=\{(x_{n})_{n\in
{\mathbb{Z}}}: x_{n}\in\mathbf{M}, x_{n+2}=ax_{n+1}+b x_{n},\forall
n\}$.
\end{rk}
It is well known (at least in the vector space case) that $\mathcal{F}(a,b)$ is a free
$\mathcal{R}$-module of rank $2$; more generally:
\begin{proposition}\label{prop1.3}
$$\mathcal{F}_{\mathbf{M}}(a,b)\cong \mathbf{M}\oplus\mathbf{M}\cong  \mathcal{F}(a,b)\otimes\mathbf{M}.$$
\end{proposition}
An explicit basis can be found for $\mathcal{F}_{\mathbf{M}}(a,b)$ (see, for example,  \cite{Niven:06} in which Lucas functions are used):
\begin{proposition}\label{prop1.4}
The sequences $\big(P_{0}^{[n]}(a,b)\big)_{n\geq\,0}$ and
$\big(P_{1}^{[n]}(a,b)\big)_{n\geq\,0}$ in
 $\mathcal{F}(a,b)$ defined by
$P_{0}^{[0]}(a,b)=1$, $P_{0}^{[1]}(a,b)=0$, respectively by
$P_{1}^{[0]}(a,b)$ $=0$, $P_{1}^{[1]}(a,b)=1$, and by $P_{i}^{[n+2]}(a,b)=aP_{i}^{[n+1]}(a,b)+bP_{i}^{[n]}(a,b)\,\,(i=0,1)$ give a canonical basis
of the $\mathcal{R}$-module $\mathcal{F}(a,b)$.
\end{proposition}
Standard operations with modules give the following:
\begin{proposition}\label{prop1.5}
 $\textbf{a})$ There is a natural $\mathcal{R}[T]$-module isomorphism:
   $$\mathcal{F}_{\mathbf{M}}(a,b)\oplus\mathcal{F}_{\mathbf{N}}(a,b)\cong\mathcal{F}_{\mathbf{M\oplus N}}(a,b)\,.$$
  $\textbf{b})$ There is a natural $\mathcal{R}^2$-module isomorphism:
  $$\mathcal{F}_{\mathbf{M}}(a,b)\oplus\mathcal{F}_{\mathbf{N}}(c,d)\cong
  \mathcal{F}_{\mathbf{M\oplus N}}\big((a,c),(b,d)\big)\,.$$
\end{proposition}

In order to describe multiplicative operations (tensor product,
symmetric power, exterior power), we introduce double Fibonacci
sequences.
\begin{defn}\label{def1.6}
  The double sequence $(x_{n,k})_{n,k\geq\,0}$, $x_{n,k}\in\mathbf{M}$ is a double Fibonacci sequence of
  type $(a,b)\otimes(c,d)\in\mathcal{R}^2\otimes\mathcal{R}^2$ if for any $n,k\geq0$ we have:
   $$x_{n+2,k}=a x_{n+1,k}+b x_{n,k}\,,$$
   $$x_{n,k+2}=c x_{n,k+1}+d x_{n,k}\,.$$
\end{defn}
As an example, let us consider the element in
$\mathcal{F}_{\mathbb{Z}}^{[2]}\big((1,1)\otimes(1,3)\big)$ with
$x_{0,0}=x_{1,0}=x_{1,1}=1$ and $x_{0,1}=0$ (we locate the terms in
the first quadrant): $$\begin{array}{ccccc}
  \vdots & \vdots & \vdots & \vdots &  \\
  3 & 7 & 10 & 17 & \cdots \\
  3 & 4 & 7 & 11 & \cdots \\
  0 & 1 & 1 & 2 & \cdots \\
  1 & 1 & 2 & 3 & \cdots
\end{array}$$
The set of double Fibonacci sequences is denoted by
$\mathcal{F}_{\mathbf{M}}^{[2]}\big((a,b)\otimes(c,d)\big)$ and it
is naturally an $\mathcal{R}[H,V]$-module \big($H$, $V$ are
horizontal and vertical shifts: $H(x_{n,k})=(x_{n+1,k})$,
respectively $V(x_{n,k})=(x_{n,k+1})$\big). If $(a,b)=(c,d)$ we use
the simplified notation $\mathcal{F}_{\mathbf{M}}^{[2]}(a,b)$.
 In \cite{Odlyzko:95} double sequences $(x_{n,k})$ given by a
 different recurrency are considered: $x_{n,k}$ depends linearly on
 the terms $\{x_{i,j}\}_{i+j<n+k}$. In our definition, $x_{n,k}$
 depends on $x_{n-1,k}$ and $x_{n-2,k}$ and also depends on $x_{n,k-1}$ and
 $x_{n,k-2}$, using two different relations. Even the existence of
 a sequence with prescribed initial four terms $x_{i,j}$, $(i,j)\in \{0,1\}^{2}$,
  is not an obvious fact. Now we present some properties and
  operations with these sequences.

In Section\,\ref{sec2} the proofs of the previous results are given. In
Section\,\ref{sec3} we generalize these results in two
directions: we consider higher order linear recurrency:
$$x_{n+d}=a_{1}x_{n+d-1}+\cdots+a_{d}x_{n},$$ and also we consider
multiple sequences:
$(x_{n_{1},n_{2},\ldots,n_{d}})_{n_{i}\geq\,0}\,.$

In the last section examples of double Fibonacci sequences are given
and also an interesting property of their diagonals is presented.
\begin{proposition}\label{prop1.7}
There is a natural isomorphism of  $\mathcal{R}[H,V]$-modules:
$$\mathcal{F}_{\mathbf{M}}(a,b)\otimes_{\mathcal{R}}\mathcal{F}_{\mathbf{N}}(c,d)
\cong\mathcal{F}_{\mathbf{M\otimes
N}}^{[2]}\big((a,b)\otimes(c,d)\big)\,.$$
\end{proposition}
\begin{cor}\label{cor1.8}
The module $\mathcal{F}^{[2]}\big((a,b)\otimes(c,d)\big)$ is a free
$\mathcal{R}$-module of rank $4$. In general,
$\mathcal{F}_{\mathbf{M\otimes N}}^{[2]}\big((a,b)\otimes(c,d)\big)$
is isomorphic to $(\mathbf{M}\otimes \mathbf{N})^{4}$.
\end{cor}
An explicit basis of $\mathcal{F}^{[2]}\big((a,b)\otimes(c,d)\big)$
is given by the four sequences
$\big(P_{i,j}^{[n,k]}(a,b)\otimes(c,d)\big)_{n,k\geq0}=\big(P_{i}^{[n]}(a,b)P_{j}^{[k]}(c,d)\big)_{n,k\geq0},$
where $(i,j)\in \{0,1\}^{2}$.\\\indent The generating function of a
double sequence $(x_{n,k})_{n,k}\geq0$ is the formal series in
$\mathcal{R}[[t,s]]\otimes\mathbf{M}\cong\mathbf{M}[[t,s]]:$
$$G(t,s)=x_{0,0}+x_{1,0}t+x_{0,1}s+\cdots+x_{n,k}t^{n}s^{k}+\cdots\,.$$

\begin{proposition}\label{prop1.9}
A Fibonacci sequence $(x_{n,k})$ of type $(a,b)\otimes(c,d)$ has a
rational generating function given by
$$G(t,s)=q(t)^{-1}r(s)^{-1}\big[x_{0,0}(1-at)(1-cs)+x_{1,0}t(1-cs)
+x_{0,1}(1-at)s+x_{1,1}ts\big]$$ where $q(t)=1-at-bt^{2}$,
$r(s)=1-cs-ds^{2}$.
\end{proposition}

\section{\large{\textbf{Proofs}}}\label{sec2}
We can write well-known results on Fibonacci sequences in the following form:
\begin{lemma}\label{lem2.1}
There are polynomials $P_{0}^{[n]},P_{1}^{[n]}\in\mathcal{R}[T,U]$
such that for any $(x_{n})_{n\geq\,0}\in
\mathcal{F}_{\mathbf{M}}(a,b):$
\begin{equation}\label{eq1}
    x_{n}=P_{0}^{[n]}(a,b)x_{0}+P_{1}^{[n]}(a,b)x_{1}
\end{equation}
for every $n\geq0$.
\end{lemma}
\begin{proof}
We define $P_{0}^{[0]}=1,\,P_{0}^{[1]}=0$ and
$P_{1}^{[0]}=0,\,P_{1}^{[1]}=1,$ and $P_{i}^{[n+2]}=a
P_{i}^{[n+1]}+bP_{i}^{[n]}\,\,(i=0,1).$ These satisfy the
equation\,(\ref{eq1}) by definition for $n=0,1$ and by induction for
$n\geq2.$
\end{proof}
\begin{rk}\label{rem2.2}
The Lemma\,\ref{lem2.1} shows that the $\mathcal{R}$-module
$\mathcal{F}(a,b)$ is free of rank $2$ with basis
$\big(P_{0}^{[n]}(a,b)\big)_{n\geq0}$,
$\big(P_{1}^{[n]}(a,b)\big)_{n\geq0}$.
\end{rk}
\begin{rk}\label{rem2.3}
If $a=r_{1}+r_{2}$, $b=-r_{1}r_{2}$ then one can describe
$P_{0}^{[n]}$ and $P_{1}^{[n]}$ in the classical way as polynomials
in $r_{1},r_{2}$:
$$\begin{array}{l}
    P_{0}^{[n]}(r_{1}+r_{2},-r_{1}r_{2})=R_{0}^{[n]}(r_{1},r_{2})=-r_{1}^{n-1}r_{2}-r_{1}^{n-2}r_{2}^{2}-\cdots
-r_{1}r_{2}^{n-1}, \\
    P_{1}^{[n]}(r_{1}+r_{2},-r_{1}r_{2})=R_{1}^{[n]}(r_{1},r_{2})=r_{1}^{n-1}+r_{1}^{n-2}r_{2}^{1}+\cdots
+r_{2}^{n-1},
  \end{array}
$$
 or as rational functions in $r_{1},r_{2}$:
\begin{equation}\label{eq2}
    R_{0}^{[n]}(r_{1},r_{2})=\frac{r_{1}^{n}r_{2}-r_{1}r_{2}^{n}}{r_{2}-r_{1}}\,,\,\,\,\,\,
    R_{1}^{[n]}(r_{1},r_{2})=\frac{r_{2}^{n}-r_{1}^{n}}{r_{2}-r_{1}}\,.
\end{equation}
\end{rk}
\begin{rk}\label{rem2.4}
The previous formulae are also correct in
$\mathcal{\widetilde{F}}_{\mathbf{M}}(a,b)$, $i.e.$ for negative
$n,$ if we extend the scalars to a suitable ring of fractions.
\end{rk}

For an arbitrary sequence $(x_{n})_{n\geq0}$ in
$\mathcal{S}_{\mathbf{M}}$ we define its \emph{generating function}
$G(t)$ as a formal series in
$\mathcal{R}[[t]]\otimes\mathbf{M}\cong\mathbf{M}[[t]]:$
$$G(t)=x_{0}+x_{1}t+x_{2}t^{2}+\cdots\,.$$
Another classical result is (see, for example, \cite{Aigner:07}):

\begin{lemma}\label{lem2.5}
The generating function of the Fibonacci sequence
$(x_{n})_{n\geq0}\in\mathcal{F}_{\mathbf{M}}(a,b)$ is the rational function
$$G(t)=\frac{(1-at)x_{0}+tx_{1}}{1-at-bt^{2}}=q(t)^{-1}\big[x_{0}+(x_{1}-ax_{0})t\big]\,,$$
where $q(t)=1-at-bt^{2}\,.$\end{lemma}
\begin{proof}{[Proposition\,\ref{prop1.4}]}
From Lemma\,\ref{lem2.1}, an arbitrary sequence
$(x_{n})_{n\geq0}\in\mathcal{F}(a,b)$ can be written as
 $(x_{n})_{n\geq0}=\big(P_{0}^{[n]}(a,b)\big)_{n\geq0}\,x_{0}+
\big(P_{1}^{[n]}(a,b)\big)_{n\geq0}\,x_{1}\,.$
\end{proof}
\begin{proof}{[Proposition\,\ref{prop1.3}]}
Define the
morphisms$$\mathcal{F}_{\mathbf{M}}(a,b)\mathop{\longrightarrow}\limits^{\varphi}\mathbf{M}\oplus
\mathbf{M}\mathop{\longrightarrow}\limits^{\psi}\mathcal{F}(a,b)\otimes\mathbf{M}
\mathop{\longrightarrow}\limits^{\eta}\mathcal{F}_{\mathbf{M}}(a,b)$$
by $$\varphi\big((x_{n})_{n\geq0}\big)=(x_{0},x_{1}),$$
$$\psi(x_{0},x_{1})=\big(P_{0}^{[n]}(a,b)\big)_{n\geq0}\otimes
x_{0}+\big(P_{1}^{[n]}(a,b)\big)_{n\geq0}\otimes x_{1},$$  and
$$\eta\big((c_{n})_{n\geq0}\otimes x\big)=(c_{n}x)_{n\geq0}.$$ It is easy to
check that $\eta\psi\varphi$, $\varphi\eta\psi$ and
$\psi\varphi\eta$ are identities, so $\varphi,\,\psi,\,\eta$ are
$\mathcal{R}$-module isomorphisms. It is also obvious that $\eta$
and $\psi\varphi$ are $\mathcal{R}[T]$-linear.
\end{proof}
\begin{proof}{[Proposition\,\ref{prop1.5}]}
There are canonical maps:
$$\Phi:\mathcal{F}_{\mathbf{M}}(a,b)\oplus\mathcal{F}_{\mathbf{N}}(a,b)
\longrightarrow\mathcal{F}_{\mathbf{M}\oplus\mathbf{N}}(a,b)$$
defined by
$$\Phi\big((x_{n})_{n\geq0},(y_{n})_{n\geq0}\big)=(x_{n},y_{n})_{n\geq0}$$ and
$$\Psi:\mathcal{F}_{\mathbf{M}}(a,b)\oplus\mathcal{F}_{\mathbf{N}}(c,d)
\longrightarrow\mathcal{F}_{\mathbf{M}\oplus\mathbf{N}}\big((a,c),(b,d)\big)$$
defined by
$$\Psi\big((x_{n})_{n\geq0},(y_{n})_{n\geq0}\big)=(x_{n},y_{n})_{n\geq0}\,.$$
Both are compatible with the shift.\end{proof}
\begin{proof}{[Proposition\,\ref{prop1.7}]}
Define the morphism of $\mathcal{R}[H,V]$-modules:
$$\Phi:\mathcal{F}_{\mathbf{M}}(a,b)\otimes\mathcal{F}_{\mathbf{N}}(c,d)
\longrightarrow\mathcal{F}_{\mathbf{M}\otimes\mathbf{N}}^{[2]}\big((a,b)\otimes(c,d)\big)$$
by$$\Phi\big((x_{n})_{n\geq0}\otimes(y_{k})_{k\geq0}\big)=(x_{n}\otimes
y_{k})_{n,k\geq0}\,.$$ The inverse morphism $\Psi$ can be
constructed using canonical bases $P_{0}^{[n]}(a,b)$, $P_{1}^{[n]}(a,b)$ of $\mathcal{F}_{\mathbf{M}}(a,b)$,  respectively $P_{0}^{[k]}(c,d)$, $P_{1}^{[k]}(c,d)$ of $\mathcal{F}_{\mathbf{N}}(c,d)$ and the corresponding basis $P_{i}^{[n]}(a,b)\otimes P_{j}^{[k]}(c,d),$ $i,j\in \{0,1\}^{2}$ of $\mathcal{F}_{\mathbf{M}}(a,b)\otimes\mathcal{F}_{\mathbf{N}}(c,d)$: if the
first four terms are given by $Z_{0,0}=\sum_{i\in I}m_{i}\otimes
n_{i}\,,$ $Z_{1,0}=\sum_{j\in J}m_{j}^{'}\otimes n_{j}^{'}\,,$
$Z_{0,1}=\sum_{h\in H}m_{h}^{''}\otimes n_{h}^{''}\,,$
$Z_{1,1}=\sum_{l\in L}m_{l}^{'''}\otimes n_{l}^{'''}\,,$ then $\Psi$
is defined by:
$$\begin{array}{rl}
    \Psi\big((Z_{n,k})_{n,k\geq0}\big)= & \sum_{i\in I}\big(P_{0}^{[n]}(a,b)m_{i}\big)_{n\geq0}\otimes
    \big(P_{0}^{[k]}(c,d)n_{i}\big)_{k\geq0}\\
    & +\sum_{j\in J}\big(P_{1}^{[n]}(a,b)m_{j}^{'}\big)_{n\geq0}\otimes
    \big(P_{0}^{[k]}(c,d)n_{j}^{'}\big)_{k\geq0}\\
     & +\sum_{h\in H}\big(P_{0}^{[n]}(a,b)m_{h}^{''}\big)_{n\geq0}\otimes
    \big(P_{1}^{[k]}(c,d)n_{h}^{''}\big)_{k\geq0}\\
    &+\sum_{l\in L}\big(P_{1}^{[n]}(a,b)m_{l}^{'''}\big)_{n\geq0}\otimes
    \big(P_{1}^{[k]}(c,d)n_{l}^{'''}\big)_{k\geq0}\,.
  \end{array}$$
\end{proof}

\begin{proof}{[Corollary\,\ref{cor1.8}]}
The proof is clear as
$\mathcal{F}^{[2]}\big((a,b)\otimes(c,d)\big)\cong
\mathcal{F}(a,b)\otimes\mathcal{F}(c,d)\cong
(\mathcal{R}\oplus\mathcal{R})\otimes(\mathcal{R}\oplus\mathcal{R})\cong\mathcal{R}^{4}$.
In general,
$\mathcal{F}_{\mathbf{M}\otimes\mathbf{N}}^{[2]}\big((a,b)\otimes(c,d)\big)\cong
\mathcal{F}_{\mathbf{M}}(a,b)\otimes\mathcal{F}_{\mathbf{N}}(c,d)\cong
(\mathbf {M}\oplus\mathbf
{M})\otimes(\mathbf{N}\oplus\mathbf{N})\cong(\mathbf
{M}\otimes\mathbf{N})^{4}$.
\end{proof}

\begin{cor}\label{cor2.6}
Using $a=r_{1}+r_{2}$, $b=-r_{1}r_{2}$, the general term $x_{n,k}$
of a sequence in
$\mathcal{F}_{\mathbf{M}\otimes\mathbf{N}}^{[2]}(a,b)$
is given by\\
$x_{n,k}=\Delta^{-2}\big[(r_{1}^{n}r_{2}-r_{1}r_{2}^{n})(r_{1}^{k}r_{2}-r_{1}r_{2}^{k})x_{0,0}
+(r_{2}^{n}-r_{1}^{n})(r_{1}^{k}r_{2}-r_{1}r_{2}^{k})x_{1,0}
\\ \indent\indent +(r_{1}^{n}r_{2}-r_{1}r_{2}^{n})(r_{2}^{k}-r_{1}^{k})x_{0,1}
+(r_{2}^{n}-r_{1}^{n})(r_{2}^{k}-r_{1}^{k})x_{1,1}\big],$\\ where
$\Delta=r_{2}-r_{1}$. This formula is correct for arbitrary integers
$n,k$ $\big($as an equality in the ring $\mathcal{R}(r_{1},r_{2})$
of rational functions$\big)$.
\end{cor}
\begin{proof}{[Proposition\,\ref{prop1.9}]}
Apply two times Lemma\,\ref{lem2.5}:
$$ \begin{array}{rl}
   G(t,s)= & \sum_{n\geq0}\big(\sum_{k\geq0}x_{n,k}s^{k}\big)t^{n} \\
   = & \sum_{n\geq0}\big[r(s)^{-1}x_{n,0}(1-cs)+r(s)^{-1}x_{n,1}s\big]t^{n} \\
   = & r(s)^{-1}\big[(1-cs)\sum_{n\geq0}x_{n,0}t^{n}+s\sum_{n\geq0}x_{n,1}t^{n}\big] \\
   = & q(t)^{-1}r(s)^{-1}\big\{(1-cs)[x_{0,0}(1-at)+x_{1,0}t]\\
   &\,+s[x_{0,1}(1-at)+x_{1,1}t]\big\}\,.
     \end{array}$$
\end{proof}

We consider also other operations with Fibonacci modules, for
example symmetric powers and exterior products (we suppose that $2$
is a unit in $\mathcal{R}$):
\begin{proposition}\label{prop2.7}
There are natural isomorphisms:
$$\emph{Symm}^{(2)}\,\mathcal{F}_{\mathbf{M}}(a,b)\cong
\big\{(x_{n,k})\in\mathcal{F}_{\mathbf{M}\otimes\mathbf{N}}^{[2]}(a,b):\,\,
x_{n,k}=x_{k,n}\,\,\,\forall\,\, k,n\geq0\big\},$$
$$\wedge^{(2)}\,\mathcal{F}_{\mathbf{M}}(a,b)\cong
\big\{(x_{n,k})\in\mathcal{F}_{\mathbf{M}\otimes\mathbf{N}}^{[2]}(a,b):\,\,
x_{n,k}=-x_{k,n}\,\,\,\forall\,\, k,n\geq0\big\}.$$
\end{proposition}
\section{\large{\textbf{Generalizations}}}\label{sec3}
First we introduce recurrency of order $d$:
\begin{defn}\label{def3.1}
  Let $\textbf{a}=(a_{1},\ldots,a_{d})$ be an element in $\mathcal{R}^{d}$.
The Fibonacci module of type $\textbf{a}$ associated to the module
$\mathbf{M}$ is the $\mathcal{R}[T]$-module:
$$\mathcal{F}_{\mathbf{M}}(\textbf{a})=
\big\{(x_{n})_{n\geq0}\in\mathcal{S}_{\mathbf{M}}:\,\,
x_{n+d}=a_{1}x_{n+d-1}+\cdots+a_{d}x_{n},\,\,\forall\,\,
n\geq0\big\}.$$
\end{defn}
Next we consider multiple Fibonacci sequences
$(x_{n_{1},\ldots,n_{p}})_{n_{i}\geq0}$ in $\mathbf{M}$:
\begin{defn}\label{def3.2}
  Let $\textbf{a}^{(1)}\in\mathcal{R}^{d_{1}},\ldots,\textbf{a}^{(p)}\in\mathcal{R}^{d_{p}}$.
The Fibonacci module of type
$(\textbf{a}^{(1)},\ldots,\textbf{a}^{(p)})$ associated to the
module $\mathbf{M}$ is the $\mathcal{R}[T_{1},\ldots,T_{p}]$-module:
$$\begin{array}{rl}
    \mathcal{F}_{\mathbf{M}}^{[p]}(\textbf{a}^{(1)},\ldots,\textbf{a}^{(p)})= & \big\{(x_{n_{1},\ldots,n_{p}})_{n_{i}
\geq0}:x_{n_{1},\ldots,n_{p}}\in{\mathbf{M}},\,\,
x_{n_{1},\ldots,n_{i}+d_{i},\ldots,n_{p}}= \\
     & \sum_{j=1}^{d_{i}}a_{j}^{(i)}x_{n_{1},\ldots,n_{i}+d_{i}-j,\ldots,n_{p}}\,\, \mbox{for}\,\, i=1,2,\ldots,p\big\}.
  \end{array}
$$
If $\textbf{a}^{(1)}=\cdots=\textbf{a}^{(p)}=\textbf{a}=(a_{1},\ldots,a_{d})$, we denote simply
$\mathcal{F}_{\mathbf{M}}^{[p]}(\textbf{a})=
\mathcal{F}_{\mathbf{M}}^{[p]}(a_{1},\ldots,a_{d})\,.$
\end{defn}
The previous results have obvious generalizations. For example:
\begin{proposition}\label{prop3.3}
$$\mathcal{F}_{\mathbf{M}}(a_{1},\ldots,a_{d})\cong \mathbf{M}^{d}\cong
 \mathcal{F}(a_{1},\ldots,a_{d})\otimes\mathbf{M}.$$
\end{proposition}
\begin{proposition}\label{prop3.4}
Fix $\textbf{\emph{a}}=(a_{1},\ldots,a_{d})\in\mathcal{R}^{d}$. The
sequences $\big(P_{i}^{[n]}(\textbf{\emph{a}})\big)_{n\geq0}$,
$i=0,\ldots,d-1$ in $\mathcal{F}(\textbf{\emph{a}})$ defined by
$P_{i}^{[j]}(\textbf{\emph{a}})=\delta_{ij}$ \emph{(}for
$j=0,\ldots,d-1 $\emph{)} give a canonical basis of
$\mathcal{F}(\textbf{\emph{a}})$.
\end{proposition}
\begin{lemma}\label{lem3.5}
The generating function of $(x_{n})_{n\geq0}$ in
$\mathcal{F}(\textbf{\emph{a}})$ is
$$G(t)=q(t)^{-1}\big[Q_{0}(t)x_{0}+Q_{1}(t)x_{1}+\cdots+Q_{d-1}(t)x_{d}\big],$$
where
$$Q_{i}(t)=t^{i}\big(1-a_{1}t-a_{2}t^{2}-\cdots-a_{d-i-1}t^{d-i-1}\big),\,\,i\in\{0,\ldots,d-1\},$$
and $q(t)=1-a_{1}t-a_{2}t^{2}-\cdots-a_{d}t^{d}\,.$
\end{lemma}
\begin{proposition}\label{prop3.6}
$$\mathcal{F}_{\mathbf{M_{1}}}(\textbf{\emph{a}}^{(1)})\otimes \cdots\otimes\mathcal{F}_{\mathbf{M_{p}}}(\textbf{\emph{a}}^{(p)})
\cong\mathcal{F}_{\mathbf{M_{1}\otimes\cdots\otimes\mathbf{M_{p}}}}^{[p]}(\textbf{\emph{a}}^{(1)},\ldots,\textbf{\emph{a}}^{(p)})
\,.$$ In particular,
$\mathcal{F}^{[p]}(\textbf{\emph{a}}^{(1)},\ldots,\textbf{\emph{a}}^{(p)})$
is free of rank $D=d_{1}d_{2}\cdots d_{p}\,.$
\end{proposition}
\begin{proposition}\label{prop3.7}
A multiple Fibonacci sequence $(x_{n_{1},\ldots ,n_{p}})$ of type
$(\textbf{\emph{a}}^{(1)},\ldots,\textbf{\emph{a}}^{(p)})$ has a
rational generating function:
$$G(t_{1},\ldots,t_{p})=q_{1}(t_{1})^{-1}\cdots q_{p}(t_{p})^{-1}
\Big[\sum_{0\leq j_{i}\leq d_{i}-1}Q_{j_{1}}^{(1)}(t_{1})\cdots
Q_{j_{p}}^{(p)}(t_{p})x_{j_{1},\ldots,j_{p}}\Big],$$ where
$q_{i}(t)=1-a_{1}^{(i)}t-\cdots-a_{d_{i}}^{(i)}t^{d_{i}}$ and
$Q_{0}^{(i)},\ldots,Q_{d_{i}-1}^{(i)}$ are like in Lemma\,\ref{lem3.5}.
\end{proposition}
For further applications in knot theory, we will use the next
specializations:
\begin{theorem}\label{th3.8}
Let $(x_{n_{1},\ldots,n_{p}})_{\geq0}$ be an element in
$\mathcal{F}_{\mathbf{M}}^{[p]}(r_{1}+r_{2},-r_{1}r_{2})$.\\
 $\textbf{a})$ The general term is given by
 $$x_{n_{1},\ldots,n_{p}}=\Delta^{-p}\sum_{0\leq j_{1},\ldots,j_{p}\leq 1}
 S_{j_{1}}^{[n_{1}]}(r_{1},r_{2})\cdots
 S_{j_{p}}^{[n_{p}]}(r_{1},r_{2})x_{j_{1},\ldots,j_{p}}\,,$$
 where $\Delta=r_{2}-r_{1}$,
 $S_{0}^{[n]}(r_{1},r_{2})=r_{1}^{n}r_{2}-r_{1}r_{2}^{n}$,
  $S_{1}^{[n]}(r_{1},r_{2})=r_{2}^{n}-r_{1}^{n}$;\\
$\textbf{b})$ the generating function of $(x_{n_{1},\ldots,n_{p}})$
is given by
 $$G(t_{1},\ldots,t_{p})=q{(t_{1})}^{-1}\cdots q{(t_{p})}^{-1}\sum_{0\leq j_{1},\ldots,j_{p}\leq 1}
 Q_{j_{1}}(t_{1})\cdots Q_{j_{p}}(t_{p})x_{j_{1},\ldots,j_{p}}\,,$$
 where $q(t)=(1-r_{1}t)(1-r_{2}t)$, $Q_{0}(t)=1-(r_{1}+r_{2})t$
 and  $Q_{1}(t)=t$.
\end{theorem}
\section{\large{\textbf{Examples}}}\label{sec4}
\begin{ex}\label{ex4.1}
Fibonacci module $\mathcal{F}_{\mathbb{Z}}^{[2]}(1,1)$: let us
analyze sequences with the first four entries
$(c_{i,j})_{(i,j)\in\{0,1\}^{2}}$ equal to $0$ or $1$. From the
sixteen possible choices there are $5$ primitive sequences:\\\\
 $B_{0}\,=\,${\small\begin{tabular}{|cc|}
  \hline
  \cline{1-2}
  0 & 0 \\
  0 & 0 \\
  \hline
\end{tabular}}\,\,,\,\,
$B_{1}\,=\,$ {\small\begin{tabular}{|cc|}
  \hline
  \cline{1-2}
  0 & 0 \\
  1 & 0 \\
  \hline
\end{tabular}}\,\,,\,\,
$B_{2}\,=\,$ {\small\begin{tabular}{|cc|}
  \hline
  \cline{1-2}
  1 & 0 \\
  0 & 1 \\
  \hline
\end{tabular}}\,\,,\,\,
$B_{3}\,=\,$ {\small\begin{tabular}{|cc|}
  \hline
  \cline{1-2}
  0 & 1  \\
  1 & 0 \\
  \hline
\end{tabular}}\,\,,\,\,
$B_{4}\,=\,$ {\small\begin{tabular}{|cc|}
  \hline
  \cline{1-2}
  1 & 0 \\
  1 & 1 \\
  \hline
\end{tabular}\,\,.}\\\\
 The others are shifts of these primitive sequences (see figure below):\\\\
 $H(B_{1})\,\,=\,$ {\small\begin{tabular}{|cc|}
  \hline
  \cline{1-2}
  0 & 0 \\
  0 & 1 \\
  \hline
\end{tabular}}\,\,,\,\,
$H^{2}(B_{1})\,=\,$ {\small\begin{tabular}{|cc|}
  \hline
  \cline{1-2}
  0 & 0 \\
  1 & 1 \\
  \hline
\end{tabular}}\,\,,\,\,
$V(B_{1})\,=\,$ {\small\begin{tabular}{|cc|}
  \hline
  \cline{1-2}
  1 & 0 \\
  0 & 0 \\
  \hline
\end{tabular}}\,\,,\,\,
$V^{2}(B_{1})\,=\,$ {\small\begin{tabular}{|cc|}
  \hline
  \cline{1-2}
  1 & 0 \\
  1 & 0 \\
  \hline
\end{tabular}}\,,
\\\\
$HV(B_{1})=$ {\small\begin{tabular}{|cc|}
  \hline
  \cline{1-2}
  0 & 1 \\
  0 & 0 \\
  \hline
\end{tabular}}\,,
$H^{2}V(B_{1})=$ {\small\begin{tabular}{|cc|}
  \hline
  \cline{1-2} 
  1 & 1 \\
  0 & 0 \\
  \hline
\end{tabular}}\,,
$HV^{2}(B_{1})=$ {\small\begin{tabular}{|cc|}
  \hline
  \cline{1-2} 
  0 & 1 \\
  0 & 1 \\
  \hline
\end{tabular}}\,,
$H^{2}V^{2}(B_{1})=$ {\small\begin{tabular}{|cc|}
  \hline
  \cline{1-2} 
  1 & 1 \\
  1 & 1 \\
  \hline
\end{tabular}}\,.\\\\

$H(B_{2})\,=\,$ {\small\begin{tabular}{|cc|}
  \hline
  \cline{1-2} 
  0 & 1 \\
  1 & 1 \\
  \hline
\end{tabular}}\,,
$V(B_{2})=$ {\small\begin{tabular}{|cc|}
  \hline
  \cline{1-2} 
  1 & 1 \\
  1 & 0 \\
  \hline
\end{tabular}}\,\,\,\,and\,
$H(B_{3})=V(B_{3})=$ {\small\begin{tabular}{|cc|}
  \hline
  \cline{1-2} 
  1 & 1 \\
  0 & 1 \\
  \hline
\end{tabular}}\,.\\\\

In fact, using the structure of $\mathbb{Z}[H,V]$-module,
$\mathcal{F}_{\mathbb{Z}}^{[2]}(1,1)$ is generated by $B_{1}$.

It is obvious that an element
$(x_{n})_{n\geq0}\in\mathcal{F}_{\mathbb{Q}}(1,1)$ can be defined by
any two terms $\{x_{p},x_{q}\}$; in the case of a double sequence
$(x_{n,k})_{n,k\geq0}\in\mathcal{F}^{[2]}_{\mathbb{Q}}(1,1)$, not any four
terms $\{x_{l,m},x_{p,q},x_{r,s},x_{u,v}\}$ can define the sequence.

\end{ex}

$$\begin{array}{ccccccccccccccccccccccccccc}
 13 &   &   &   &   &   &   &   &    &   &   &   &   & 21 &   &   &   &   &   &   &    & \\
  8 & 0 &   &   &   &   &   &   &    &   &   &   &   & 13 & 8 &   &   &   &   &   &    & \\
  5 & 0 & 5 &   &   &   &   &   &    &   &   &   &   & 8  & 5 & 13&   &   &   &   &    & \\
  3 & 0 & 3 & 3 &   &   &   &   &    &   &   &   &   & 5  & 3 & 8 & 11&   &   &   &    & \\
  2 & 0 & 2 & 2 & 4 &   &   &   &    &   &   &   &   & 3  & 2 & 5 & 7 & 12&   &   &    & \\
  1 & 0 & 1 & 1 & 2 & 3 &   &   &    &   &   &   &   & 2  & 1 & 3 & 4 & 7 & 11&   &    & \\
  1 & 0 & 1 & 1 & 2 & 3 & 5 &   &    &   &   &   &   & 1  & 1 & 2 & 3 & 5 & 8 &13 &    & \\
  0 & 0 & 0 & 0 & 0 & 0 & 0 & 0 &    &   &   &   &   & 1  & 0 & 1 & 1 & 2 & 3 & 5 & 8  & \\
  1 & 0 & 1 & 1 & 2 & 3 & 5 & 8 & 13 &   &   &   &   & 0  & 1 & 1 & 2 & 3 & 5 & 8 & 13 & 21
\end{array}$$

\bigskip
A curious property of these sequences is the alternating
monotonicity along the lines parallel to the secondary diagonal:
$$x_{n+2,k}\geq x_{n+1,k+1}\leq x_{n,k+2}$$ or $$x_{n+2,k}\leq x_{n+1,k+1}\geq x_{n,k+2}\,.$$
In general we do not have this strong alternating property (look at the
sequence given by $x_{0,0}=x_{1,0}=3,x_{0,1}=2, x_{1,1}=0$: the 4th
diagonal is $(7,3,2,9)$ ). In general we have only a "weak alternating
property": $$x_{n+2,k+1}\geq x_{n+1,k+2}\,\, \mbox{if and only if}
\,\,x_{n+3,k}\leq x_{n,k+3}$$ (see the next corollary).

The general statement explaining these two facts is given by:
\begin{proposition}{\emph{(diagonal property)}}\label{prop4.2}
If $a^{2}d=bc^{2}$, any four diagonal consecutive terms of the
sequence
$(x_{n,k})_{n,k\geq0}\in\mathcal{F}_{\mathbf{M}}^{[2]}\big((a,b)\otimes(c,d)\big)$
satisfy the relation:
$$abx_{n,k+3}+(a^{2}+b)cx_{n+1,k+2}=a(c^{2}+d)x_{n+2,k+1}+cdx_{n+3,k}\,.$$
\end{proposition}
\begin{proof}
Express the terms as combinations of $x_{n,k},$
$x_{n+1,k},$ $x_{n,k+1}$ and $x_{n+1,k+1}$.
\end{proof}
\begin{cor}\label{cor4.3}
Four diagonal consecutive terms in
$(x_{n,k})_{n,k\geq0}\in\mathcal{F}^{[2]}_{\mathbb{Z}}(1,1)$ satisfy
$$x_{n,k+3}-x_{n+3,k}=2(x_{n+2,k+1}-x_{n+1,k+2})\,.$$
\end{cor}
\flushleft{\textbf{Acknowledgment.}}
I would like to thank the referee for many comments and improvements of the first version of the paper.

\medskip
\textbf{2000 Mathematics Subject Classification: 05A15, 11B39}


\begin{thebibliography}{10}
\bibitem{Aigner:07} M. Aigner, \emph{A Course in Enumeration}, Springer 2007.
\bibitem{Niven:06} I. Niven, H. Montgomery, \emph{An Introduction to the Theory of Numbers},
John Wiley and Sons, 2006.
\bibitem{Odlyzko:95} A. M. Odlyzko, \emph{Asymptotic Enumeration Methods}, in R. L. Graham,
M. Gr\"{o}tschel, L. Lov\'{a}sz: Handbook of Combinatorics, Vol.II,
pp.1063-1230, Elsevier 1995.
\end{thebibliography}
\end{document}